# FAST LEARNING RATES FOR PLUG-IN CLASSIFIERS

BY JEAN-YVES AUDIBERT AND ALEXANDRE B. TSYBAKOV

*École Nationale des Ponts et Chaussées and Université Paris VI*

It has been recently shown that, under the margin (or low noise) assumption, there exist classifiers attaining fast rates of convergence of the excess Bayes risk, that is, rates faster than $n^{-1/2}$. The work on this subject has suggested the following two conjectures: (i) the best achievable fast rate is of the order $n^{-1}$, and (ii) the plug-in classifiers generally converge more slowly than the classifiers based on empirical risk minimization. We show that both conjectures are not correct. In particular, we construct plug-in classifiers that can achieve not only fast, but also *super-fast* rates, that is, rates faster than $n^{-1}$. We establish minimax lower bounds showing that the obtained rates cannot be improved.

**1. Introduction.** Let $(X, Y)$ be a random couple taking values in $\mathcal{Z} \triangleq \mathbf{R}^d \times \{0, 1\}$ with joint distribution $P$. We regard $X \in \mathbf{R}^d$ as a vector of features corresponding to an object and $Y \in \{0, 1\}$ as a label indicating that the object belongs to one of two classes. Consider the sample $(X_1, Y_1), \ldots, (X_n, Y_n)$, where $(X_i, Y_i)$ are independent copies of $(X, Y)$. We denote by $P^{\otimes n}$ the product probability measure according to which the sample is distributed, and by $P_X$ the marginal distribution of $X$.

The goal of a classification procedure is to predict the label $Y$ given the value of $X$, that is, to provide a decision rule $f : \mathbf{R}^d \to \{0, 1\}$ which belongs to the set $\mathcal{F}$ of all Borel functions defined on $\mathbf{R}^d$ and taking values in $\{0, 1\}$. The performance of a decision rule $f$ is measured by the misclassification error

$$R(f) \triangleq P(Y \neq f(X)).$$

The Bayes decision rule is a minimizer of the risk $R(f)$ over all decision rules $f \in \mathcal{F}$, and one of such minimizers has the form $f^*(X) = \mathbb{1}_{\{\eta(X) \geq 1/2\}}$, where $\mathbb{1}_{\{\cdot\}}$ denotes the indicator function and $\eta(X) \stackrel{\triangle}{=} P(Y = 1|X)$ is the regression function of $Y$ on $X$ [here $P(dY|X)$ is a regular conditional probability, which we will use in the following without further mention].

An empirical decision rule (a classifier) is a random mapping $\hat{f}_n : \mathcal{Z}^n \to \mathcal{F}$ measurable w.r.t. the sample. Its accuracy can be characterized by the excess risk,

$$(1.1) \qquad \mathcal{E}(\hat{f}_n) = \mathbf{E} R(\hat{f}_n) - R(f^*) = \mathbf{E}(|2\eta(X) - 1|\mathbb{1}_{\{\hat{f}_n(X) \neq f^*(X)\}})$$

where $\mathbf{E}$ denotes expectation. A key problem in classification is to construct classifiers with small excess risk (cf. [8, 24]). Optimal classifiers can be defined as those having the best possible rate of convergence of $\mathcal{E}(\hat{f}_n)$ to 0, as $n \to \infty$. Of course, this rate, and thus the optimal classifier, depend on the assumptions on the joint distribution of $(X, Y)$. A standard way to define optimal classifiers is to introduce a class of joint distributions of $(X, Y)$ and to declare $\hat{f}_n$ optimal if it achieves the best rate of convergence in a minimax sense on this class.

Two types of assumptions on the joint distribution of $(X, Y)$ are commonly used: complexity assumptions and margin assumptions.

*Complexity assumptions* are stated in two possible ways. The first of them is to suppose that the regression function $\eta$ is smooth enough or, more generally, belongs to a class of functions $\Sigma$ having a suitably bounded $\varepsilon$-entropy. This is called a *complexity assumption on the regression function* (CAR). Most commonly it is of the following form.

ASSUMPTION (CAR). The regression function $\eta$ belongs to the class $\Sigma$ of functions on $\mathbf{R}^d$ such that

$$\mathcal{H}(\varepsilon, \Sigma, L_p) \leq A_* \varepsilon^{-\rho} \qquad \forall \varepsilon > 0,$$

with some constants $\rho > 0$, $A_* > 0$. Here $\mathcal{H}(\varepsilon, \Sigma, L_p)$ denotes the $\varepsilon$-entropy of the set $\Sigma$ w.r.t. an $L_p$ norm with some $1 \leq p \leq \infty$.

Recall that the metric entropy $\mathcal{H}(\varepsilon, \Sigma, L_p)$ is the logarithm of the minimum number of $L_p$-balls of radius $\varepsilon$ covering the set $\Sigma$ [10].

At this stage of discussion we do not identify precisely the value of $p$ for the $L_p$ norm in Assumption (CAR), or the measure with respect to which this norm is defined. Examples will be given later. If $\Sigma$ is a class of smooth functions with smoothness parameter $\beta$ on a compact in $\mathbf{R}^d$, for example, a Hölder class, as described below, a typical value of $\rho$ in Assumption (CAR) is $\rho = d/\beta$.

Assumption (CAR) is well adapted for the study of *plug-in rules*, that is, of the classifiers having the form

$$(1.2) \qquad \hat{f}_n^{\text{PI}}(X) = \mathbb{1}_{\{\hat{\eta}_n(X) \geq 1/2\}},$$

where $\hat{\eta}_n$ is a nonparametric estimator of the function $\eta$. Indeed, Assumption (CAR) typically reads as a smoothness assumption on $\eta$, implying that a good nonparametric estimator (kernel, local polynomial, orthogonal series or other) $\hat{\eta}_n$ converges with some rate to the regression function $\eta$, as $n \to \infty$. In turn, closeness of $\hat{\eta}_n$ to $\eta$ implies closeness of $\hat{f}_n$ to $f$: for any plug-in classifier $\hat{f}_n^{\mathrm{PI}}$ we have

$$(1.3) \qquad \mathbf{E} R(\hat{f}_n^{\mathrm{PI}}) - R(f^*) \leq 2\mathbf{E} \int |\hat{\eta}_n(x) - \eta(x)| P_X(dx)$$

(cf. [8], Theorem 2.2). For various types of estimators $\hat{\eta}_n$ and under rather general assumptions it can be shown that, if Assumption (CAR) holds, the RHS of (1.3) is uniformly of the order $n^{-1/(2+\rho)}$, and thus

$$(1.4) \qquad \sup_{P:\eta \in \Sigma} \mathcal{E}(\hat{f}_n^{\mathrm{PI}}) = O(n^{-1/(2+\rho)}), \qquad n \to \infty$$

(cf. [26]). In particular, if $\rho = d/\beta$ (which corresponds to a class of smooth functions with smoothness parameter $\beta$), we get

$$(1.5) \qquad \sup_{P:\eta \in \Sigma} \mathcal{E}(\hat{f}_n^{\mathrm{PI}}) = O(n^{-\beta/(2\beta+d)}), \qquad n \to \infty.$$

Note that (1.5) can be easily deduced from (1.3) and standard results on the $L_1$ or $L_2$ convergence rates of usual nonparametric regression estimators on $\beta$-smoothness classes $\Sigma$. The rates in (1.4), (1.5) are quite slow, always slower than $n^{-1/2}$. In (1.5) they deteriorate dramatically as the dimension $d$ increases. Moreover, Yang [26] showed that, under general assumptions, the bound (1.5) cannot be improved in a minimax sense. These results raised some pessimism about the plug-in rules.

The second way to describe complexity is to introduce a structure on the class of possible decision sets $G^* = \{x : f^*(x) = 1\} = \{x : \eta(x) \geq 1/2\}$ rather than on that of regression functions $\eta$. A standard *complexity assumption on the decision set* (CAD) is the following.

ASSUMPTION (CAD). The decision set $G^*$ belongs to a class $\mathcal{G}$ of subsets of $\mathbf{R}^d$ such that

$$\mathcal{H}(\varepsilon, \mathcal{G}, d_\triangle) \leq A_* \varepsilon^{-\rho} \qquad \forall \varepsilon > 0,$$

with some constants $\rho > 0$, $A_* > 0$. Here $\mathcal{H}(\varepsilon, \mathcal{G}, d_\triangle)$ denotes the $\varepsilon$-entropy of the class $\mathcal{G}$ w.r.t. the measure of symmetric difference pseudo-distance between sets defined by $d_\triangle(G, G') = P_X(G \triangle G')$ for two measurable subsets $G$ and $G'$ in $\mathbf{R}^d$.

The parameter $\rho$ in Assumption (CAD) typically characterizes the smoothness of the boundary of $G^*$ (cf. [20]). Note that, in general, there is no connection between Assumptions (CAR) and (CAD). Indeed, the fact that $G^*$ has a smooth boundary does not imply that $\eta$ is smooth, and vice versa. In Assumption (CAD), the values of $\rho$ closer to 0 correspond to smoother boundaries (less complex sets $G^*$). As a limit case when $\rho \to 0$ one can consider the Vapnik–Chervonenkis classes (VC-classes) for which the $\varepsilon$-entropy is logarithmic in $1/\varepsilon$.

Assumption (CAD) is suited for the study of empirical risk minimization (ERM) type classifiers introduced by Vapnik and Chervonenkis [25]; see also [8, 24]. As shown in [20], for every $0 < \rho < 1$ there exist ERM classifiers $\hat{f}_n^{\mathrm{ERM}}$ such that, under Assumption (CAD),

$$(1.6) \qquad \sup_{P\,:\,G^* \in \mathcal{G}} \mathcal{E}(\hat{f}_n^{\mathrm{ERM}}) = O(n^{-1/2}), \qquad n \to \infty.$$

The rate of convergence in (1.6) is better than that for plug-in rules [see (1.4) and (1.5)] and it does not depend on $\rho$ (resp., on the dimension $d$). Note that the comparison between (1.6) and (1.4) and (1.5) is not quite legitimate, because classes of joint distributions $P$ of $(X, Y)$ satisfying Assumption (CAR) are different from those satisfying Assumption (CAD). Nevertheless, such a comparison has been often interpreted as an argument in disfavor of the plug-in rules. Indeed, Yang's lower bound shows that the $n^{-1/2}$ rate cannot be attained under Assumption (CAR) suited for the plug-in rules. Recently, advantages of the ERM type classifiers, including penalized ERM methods, have been further confirmed by the fact that, under the margin (or low noise) assumption, they can attain *fast rates of convergence*, that is, rates that are faster than $n^{-1/2}$ [1, 11, 14, 15, 20, 22].

The *margin assumption* (or low noise assumption) is stated as follows.

ASSUMPTION (MA). There exist constants $C_0 > 0$ and $\alpha \geq 0$ such that

$$(1.7) \qquad P_X(0 < |\eta(X) - 1/2| \leq t) \leq C_0 t^\alpha \qquad \forall t > 0.$$

The case $\alpha = 0$ is trivial (no assumption) and is included for notational convenience. The other extreme case $\alpha = \infty$ is most advantageous for classification: the regression function $\eta$ is bounded away from $1/2$. Assumption (MA) provides a useful characterization of the behavior of the regression function $\eta$ in the vicinity of the level $\eta = 1/2$, which turns out to be crucial for convergence of classifiers. Note that the margin assumption does not affect the complexity of the class of regression functions, but it affects the rate of convergence of the excess risk due to its structure. This can be seen from the following simple argument which underlies our results. For any $\delta > 0$ from (1.1) and Assumption (MA) we get

$$\begin{aligned}
\mathbf{E} R(\hat{f}_n^{\mathrm{PI}}) - R(f^*) &\leq 2\delta P_X(0 < |\eta(X) - 1/2| \leq \delta) \\
(1.8) &\qquad + \mathbf{E}(|2\eta(X) - 1| \mathbb{1}_{\{\hat{f}_n^{\mathrm{PI}}(X) \neq f^*(X)\}} \mathbb{1}_{\{|\eta(X) - 1/2| > \delta\}}) \\
&\leq 2 C_0 \delta^{1+\alpha} + 2\mathbf{E}(|\hat{\eta}_n(X) - \eta(X)| \mathbb{1}_{\{|\hat{\eta}_n(X) - \eta(X)| > \delta\}}),
\end{aligned}$$

where in the last inequality we have used the fact that $|\eta(X) - 1/2| \leq |\hat{\eta}_n(X) - \eta(X)|$ on the set $\{X : \hat{f}_n^{\text{PI}}(X) \neq f^*(X)\}$. Thus, the excess risk is decomposed into two terms: $2C_0\delta^{1+\alpha}$ which is determined by the margin assumption and reflects the behavior of $\eta$ near the decision boundary, and the second term which characterizes the regression estimation error. Optimal convergence is essentially obtained by choosing the $\delta$ which balances the two terms. Fast rates are possible because the second term in (1.8) decreases exponentially in $\delta$ for several types of regression estimators $\hat{\eta}_n$.

For more discussion of the margin assumption see [20] and the survey [7]. The main point is that, under Assumption (MA), fast classification rates up to $n^{-1}$ are achievable. In particular, for every $0 < \rho < 1$ and $\alpha > 0$ there exist ERM type classifiers $\hat{f}_n^{\text{ERM}}$ such that

$$(1.9) \qquad \sup_{P:\,(\text{CAD}),(\text{MA})} \mathcal{E}(\hat{f}_n^{\text{ERM}}) = O(n^{-(1+\alpha)/(2+\alpha+\alpha\rho)}), \qquad n \to \infty,$$

where $\sup_{P:\,(\text{CAD}),(\text{MA})}$ denotes the supremum over all joint distributions $P$ of $(X,Y)$ satisfying Assumptions (CAD) and (MA). The RHS of (1.9) can be arbitrarily close to $O(n^{-1})$ for large $\alpha$ and small $\rho$. Result (1.9) for direct ERM classifiers on $\varepsilon$-nets is proved by Tsybakov [20], and for some other ERM type classifiers by Tsybakov and van de Geer [22], Koltchinskii [11] and Audibert [1] [in some of these papers the rate of convergence (1.9) is obtained with an extra log-factor].

Comparison of (1.6) and (1.9) with (1.4) seems to support the conjecture that the plug-in classifiers are inferior to the ERM type ones. The main message of the present paper is to disprove this. We will show that there exist plug-in rules converging with fast rates, and even with *super-fast rates*, that is, faster than $n^{-1}$ under the margin Assumption (MA). The basic idea of the proof is to use arguments similar to (1.8) combined with exponential inequalities for the regression estimator $\hat{\eta}_n$ (see Section 3 below) or the convergence results in the $L_\infty$ norm (see Section 5), rather than the usual $L_1$ or $L_2$ norm convergence of $\hat{\eta}_n$ as previously described [cf. (1.3)]. On the other hand, the super-fast rates are not attainable for ERM type rules or, more precisely, under Assumption (CAD), which serves for the study of ERM type rules. In fact, the lower bound of [15] shows that the rates cannot be faster than $(\log n)/n$ even for smaller classes than those satisfying (CAD).

It is important to note that our results on fast rates cover more general settings than just classification with plug-in rules. These are rather results about *classification in the regression complexity context under the margin assumption*. In particular, we establish minimax lower bounds valid for all classifiers, and we construct a "hybrid" plug-in/ERM procedure (i.e., a procedure performing ERM on a set of plug-in rules coming from an appropriate grid on the set of regression functions) that achieves optimality.

Thus, the point is mainly not about the type of procedure (plug-in or ERM) but about the type of complexity assumption [on the regression function (CAR) or on the decision set (CAD)] that should be natural to impose. Assumption (CAR) on the regression function arises in a natural way in the analysis of several practical procedures of plug-in type, such as various versions of boosting or SVM (cf. [3, 5, 6, 17, 19]). These procedures are now being intensively studied, but, to our knowledge, only suboptimal rates of convergence have been proved in the regression complexity context under the margin assumption. The results in Section 4 (see also Section 5) establish the optimal rates of classification under Assumption (CAR) and show that they are attained for a "hybrid" plug-in/ERM procedure. Expectedly, the same rates should be achievable for other plug-in type procedures, such as boosting.

**2. Notation and definitions.** In this section we introduce some notation, definitions and basic facts that will be used in the paper.

We denote by $C, C_1, C_2, \ldots$ positive constants whose values may differ from line to line. The symbols $\mathbf{P}$ and $\mathbf{E}$ stand for generic probability and expectation, and $E_X$ is the expectation w.r.t. the marginal distribution $P_X$. We denote by $\mathcal{B}(x, r)$ the closed Euclidean ball in $\mathbf{R}^d$ centered at $x \in \mathbf{R}^d$ and of radius $r > 0$.

For any multi-index $s = (s_1, \ldots, s_d) \in \mathbf{N}^d$ and any $x = (x_1, \ldots, x_d) \in \mathbf{R}^d$, we define $|s| = \sum_{i=1}^d s_i$, $s! = s_1! \cdots s_d!$, $x^s = x_1^{s_1} \cdots x_d^{s_d}$ and $\|x\| \triangleq (x_1^2 + \cdots + x_d^2)^{1/2}$. Let $D^s$ denote the differential operator $D^s \triangleq \frac{\partial^{s_1 + \cdots + s_d}}{\partial x_1^{s_1} \cdots \partial x_d^{s_d}}$.

Let $\beta > 0$. Denote by $\lfloor \beta \rfloor$ the maximal integer that is strictly less than $\beta$. For any $x \in \mathbf{R}^d$ and any $\lfloor \beta \rfloor$-times continuously differentiable real-valued function $g$ on $\mathbf{R}^d$, we denote by $g_x$ its Taylor polynomial of degree $\lfloor \beta \rfloor$ at point $x$,

$$g_x(x') \triangleq \sum_{|s| \leq \lfloor \beta \rfloor} \frac{(x' - x)^s}{s!} D^s g(x).$$

Let $L > 0$. The $(\beta, L, \mathbf{R}^d)$-*Hölder class* of functions, denoted $\Sigma(\beta, L, \mathbf{R}^d)$, is defined as the set of functions $g : \mathbf{R}^d \to \mathbf{R}$ that are $\lfloor \beta \rfloor$ times continuously differentiable and satisfy, for any $x, x' \in \mathbf{R}^d$, the inequality

$$|g(x') - g_x(x')| \leq L \|x - x'\|^\beta.$$

Fix some constants $c_0, r_0 > 0$. We will say that a Lebesgue measurable set $A \subset \mathbf{R}^d$ is $(c_0, r_0)$-*regular* if

(2.1) $\qquad \lambda[A \cap \mathcal{B}(x, r)] \geq c_0 \lambda[\mathcal{B}(x, r)] \qquad \forall \, 0 < r \leq r_0, \forall x \in A,$

where $\lambda[S]$ stands for the Lebesgue measure of $S \subset \mathbf{R}^d$. To illustrate this definition, consider the following example. Let $d \geq 2$. Then the set $A = \{x =$ $(x_1, \ldots, x_d) \in \mathbf{R}^d : \sum_{j=1}^d |x_j|^q \leq 1\}$ is $(c_0, r_0)$-regular with some $c_0, r_0 > 0$ for $q \geq 1$, and there are no $c_0, r_0 > 0$ such that $A$ is $(c_0, r_0)$-regular for $0 < q < 1$.

Introduce now two assumptions on the marginal distribution $P_X$ that will be used in the sequel.

DEFINITION 2.1. Fix $0 < c_0, r_0, \mu_{\max} < \infty$ and a compact $\mathcal{C} \subset \mathbf{R}^d$. We say that the *mild density assumption* is satisfied if the marginal distribution $P_X$ is supported on a compact $(c_0, r_0)$-regular set $A \subseteq \mathcal{C}$ and has a uniformly bounded density $\mu$ w.r.t. the Lebesgue measure: $\mu(x) \leq \mu_{\max}$, $\forall x \in A$.

DEFINITION 2.2. Fix some constants $c_0, r_0 > 0$ and $0 < \mu_{\min} < \mu_{\max} < \infty$ and a compact $\mathcal{C} \subset \mathbf{R}^d$. We say that the *strong density assumption* is satisfied if the marginal distribution $P_X$ is supported on a compact $(c_0, r_0)$-regular set $A \subseteq \mathcal{C}$ and has a density $\mu$ w.r.t. the Lebesgue measure bounded away from zero and infinity on $A$:

$$\mu_{\min} \leq \mu(x) \leq \mu_{\max} \quad \text{for } x \in A \quad \text{and} \quad \mu(x) = 0 \quad \text{otherwise.}$$

We finally recall some notions related to local polynomial estimators.

DEFINITION 2.3. For $h > 0$, $x \in \mathbf{R}^d$, for an integer $l \geq 0$ and a function $K : \mathbf{R}^d \to \mathbf{R}_+$, denote by $\hat{\theta}_x$ a polynomial on $\mathbf{R}^d$ of degree $l$ which minimizes

$$(2.2) \qquad \sum_{i=1}^n [Y_i - \hat{\theta}_x(X_i - x)]^2 K\left(\frac{X_i - x}{h}\right).$$

The *local polynomial estimator* $\hat{\eta}_n^{\mathrm{LP}}(x)$ *of order* $l$, or LP($l$) estimator, of the value $\eta(x)$ of the regression function at point $x$ is defined by $\hat{\eta}_n^{\mathrm{LP}}(x) \triangleq \hat{\theta}_x(0)$ if $\hat{\theta}_x$ is the unique minimizer of (2.2) and $\hat{\eta}_n^{\mathrm{LP}}(x) \triangleq 0$ otherwise. The value $h$ is called the bandwidth and the function $K$ is called the kernel of the LP($l$) estimator.

Let $T_s$ denote the coefficients of $\hat{\theta}_x$ indexed by the multi-index $s \in \mathbf{N}^d$, $\hat{\theta}_x(u) = \sum_{|s| \leq l} T_s u^s$. Introduce the vectors $T \triangleq (T_s)_{|s| \leq l}$, $V \triangleq (V_s)_{|s| \leq l}$, where

$$(2.3) \qquad V_s \triangleq \sum_{i=1}^n Y_i (X_i - x)^s K\left(\frac{X_i - x}{h}\right),$$

$U(u) \triangleq (u^s)_{|s| \leq l}$ and the matrix $Q \triangleq (Q_{s_1, s_2})_{|s_1|, |s_2| \leq l}$, where

$$(2.4) \qquad Q_{s_1, s_2} \triangleq \sum_{i=1}^n (X_i - x)^{s_1 + s_2} K\left(\frac{X_i - x}{h}\right).$$

The following result is straightforward (cf. Section 1.7 in [21] where the case $d = 1$ is considered).

PROPOSITION 2.1. *If the matrix $Q$ is positive definite, there exists a unique polynomial on $\mathbf{R}^d$ of degree $l$ minimizing (2.2). Its vector of coefficients is given by $T = Q^{-1}V$ and the corresponding $\mathrm{LP}(l)$ regression function estimator has the form*

$$\hat{\eta}_n^{\mathrm{LP}}(x) = U^T(0)Q^{-1}V = \sum_{i=1}^n Y_i K\left(\frac{X_i - x}{h}\right) U^T(0)Q^{-1}U(X_i - x).$$

**3. Fast rates for plug-in rules under the strong density assumption.** We first state a general result showing how the rates of convergence of plug-in classifiers can be deduced from exponential inequalities for the corresponding regression estimators.

LEMMA 3.1. *Let $\hat{\eta}_n$ be an estimator of the regression function $\eta$ and $\mathcal{P}$ a set of probability distributions on $\mathcal{Z}$ such that for some constants $C_1 > 0$, $C_2 > 0$, for some positive sequence $a_n$, for $n \geq 1$ and any $\delta > 0$, and for almost all $x$ w.r.t. $P_X$, we have*

$$(3.1) \qquad \sup_{P \in \mathcal{P}} P^{\otimes n}(|\hat{\eta}_n(x) - \eta(x)| \geq \delta) \leq C_1 \exp(-C_2 a_n \delta^2).$$

*Consider the plug-in classifier $\hat{f}_n = \mathbb{1}_{\{\hat{\eta}_n \geq 1/2\}}$. If all the distributions $P \in \mathcal{P}$ satisfy the margin Assumption (MA), we have*

$$\sup_{P \in \mathcal{P}} \{\mathbf{E} R(\hat{f}_n) - R(f^*)\} \leq C a_n^{-(1+\alpha)/2}$$

*for $n \geq 1$ with some constant $C > 0$ depending only on $\alpha$, $C_0$, $C_1$ and $C_2$.*

PROOF. Consider the sets $A_j \subset \mathbf{R}^d, j = 1, 2, \ldots$, defined as

$$A_0 \triangleq \{x \in \mathbf{R}^d : 0 < |\eta(x) - \tfrac{1}{2}| \leq \delta\},$$
$$A_j \triangleq \{x \in \mathbf{R}^d : 2^{j-1}\delta < |\eta(x) - \tfrac{1}{2}| \leq 2^j \delta\} \qquad \text{for } j \geq 1.$$

For any $\delta > 0$, we may write

$$\mathbf{E} R(\hat{f}_n) - R(f^*) = \mathbf{E}(|2\eta(X) - 1|\mathbb{1}_{\{\hat{f}_n(X) \neq f^*(X)\}})$$

$$= \sum_{j=0}^\infty \mathbf{E}(|2\eta(X) - 1|\mathbb{1}_{\{\hat{f}_n(X) \neq f^*(X)\}}\mathbb{1}_{\{X \in A_j\}})$$

(3.2)

$$\leq 2\delta P_X(0 < |\eta(X) - \tfrac{1}{2}| \leq \delta)$$
$$+ \sum_{j \geq 1} \mathbf{E}(|2\eta(X) - 1|\mathbb{1}_{\{\hat{f}_n(X) \neq f^*(X)\}}\mathbb{1}_{\{X \in A_j\}}).$$

On the event $\{\hat{f}_n \neq f^*\}$ we have $|\eta - \frac{1}{2}| \leq |\hat{\eta}_n - \eta|$. So, for any $j \geq 1$, we get

$$\mathbf{E}(|2\eta(X) - 1|\mathbb{1}_{\{\hat{f}_n(X) \neq f^*(X)\}}\mathbb{1}_{\{X \in A_j\}})$$
$$\leq 2^{j+1}\delta \mathbf{E}[\mathbb{1}_{\{|\hat{\eta}_n(X) - \eta(X)| \geq 2^{j-1}\delta\}}\mathbb{1}_{\{0 < |\eta(X) - 1/2| \leq 2^j\delta\}}]$$
$$\leq 2^{j+1}\delta E_X[P^{\otimes n}(|\hat{\eta}_n(X) - \eta(X)| \geq 2^{j-1}\delta)\mathbb{1}_{\{0 < |\eta(X) - 1/2| \leq 2^j\delta\}}]$$
$$\leq C_1 2^{j+1}\delta \exp(-C_2 a_n (2^{j-1}\delta)^2) P_X(0 < |\eta(X) - \tfrac{1}{2}| \leq 2^j\delta)$$
$$\leq 2C_1 C_0 2^{j(1+\alpha)}\delta^{1+\alpha} \exp(-C_2 a_n (2^{j-1}\delta)^2),$$

where in the last inequality we have used Assumption (MA). Now, from inequality (3.2), taking $\delta = a_n^{-1/2}$ and using Assumption (MA) to bound the first term of the right-hand side of (3.2), we get

$$\mathbf{E}R(\hat{f}_n) - R(f^*) \leq 2C_0 a_n^{-(1+\alpha)/2} + C a_n^{-(1+\alpha)/2} \sum_{j \geq 2} 2^{j(1+\alpha)} \exp(-C_2 2^{2j-2})$$
$$\leq C a_n^{-(1+\alpha)/2}. \qquad \square$$

Inequality (3.1) is crucial to obtain the above result. This inequality holds for various types of estimators and various sets of probability distributions $\mathcal{P}$. Here we focus on the standard case where $\eta$ belongs to the Hölder class $\Sigma(\beta, L, \mathbf{R}^d)$ and the marginal law of $X$ satisfies the strong density assumption. We are going to show that in this case there exist estimators satisfying inequality (3.1) with $a_n = n^{2\beta/(2\beta+d)}$. These can be, for example, locally polynomial estimators. Specifically, assume from now on that $K$ is a kernel satisfying

(3.3) $$\exists c > 0: \quad K(x) \geq c\mathbb{1}_{\{\|x\| \leq c\}} \quad \forall x \in \mathbf{R}^d,$$

(3.4) $$\int_{\mathbf{R}^d} K(u)\,du = 1,$$

(3.5) $$\int_{\mathbf{R}^d} (1 + \|u\|^{4\beta}) K^2(u)\,du < \infty,$$

(3.6) $$\sup_{u \in \mathbf{R}^d} (1 + \|u\|^{2\beta}) K(u) < \infty.$$

Let $h > 0$, and consider the matrix $\bar{B} \triangleq (\bar{B}_{s_1,s_2})_{|s_1|,|s_2| \leq \lfloor\beta\rfloor}$, where $\bar{B}_{s_1,s_2} = \frac{1}{nh^d}\sum_{i=1}^n (\frac{X_i - x}{h})^{s_1+s_2} K(\frac{X_i - x}{h})$. Define the regression function estimator $\hat{\eta}_n^*$ as follows. If the smallest eigenvalue of the matrix $\bar{B}$ is greater than $(\log n)^{-1}$ we set $\hat{\eta}_n^*(x)$ equal to the projection of $\hat{\eta}_n^{\mathrm{LP}}(x)$ on the interval $[0,1]$, where $\hat{\eta}_n^{\mathrm{LP}}(x)$ is the $\mathrm{LP}(\lfloor\beta\rfloor)$ estimator with bandwidth $h > 0$ and kernel $K$ satisfying (3.3)–(3.6). If the smallest eigenvalue of $\bar{B}$ is less than $(\log n)^{-1}$ we set $\hat{\eta}_n^*(x) = 0$.

THEOREM 3.2. *Let $\mathcal{P}$ be a class of probability distributions $P$ on $\mathcal{Z}$ such that the regression function $\eta$ belongs to the Hölder class $\Sigma(\beta, L, \mathbf{R}^d)$ and the marginal law of $X$ satisfies the strong density assumption. Then there exist constants $C_1, C_2, C_3 > 0$ such that for any $0 < h \leq r_0/c$, any $C_3 h^\beta < \delta$ and any $n \geq 1$ the estimator $\hat{\eta}_n^*$ satisfies*

$$\sup_{P \in \mathcal{P}} P^{\otimes n}(|\hat{\eta}_n^*(x) - \eta(x)| \geq \delta) \leq C_1 \exp(-C_2 n h^d \delta^2) \tag{3.7}$$

*for almost all $x$ w.r.t. $P_X$. As a consequence, there exist $C_1, C_2 > 0$ such that for $h = n^{-1/(2\beta+d)}$ and any $\delta > 0$, $n \geq 1$ we have*

$$\sup_{P \in \mathcal{P}} P^{\otimes n}(|\hat{\eta}_n^*(x) - \eta(x)| \geq \delta) \leq C_1 \exp(-C_2 n^{2\beta/(2\beta+d)} \delta^2) \tag{3.8}$$

*for almost all $x$ w.r.t. $P_X$. The constants $C_1, C_2, C_3$ depend only on $\beta$, $d$, $L$, $c_0$, $r_0$, $\mu_{\min}, \mu_{\max}$, and on the kernel $K$.*

The proof is given in Section 6.1.

REMARK 3.1. We have chosen here the LP estimators of $\eta$, because for them the exponential inequality (3.1) holds without additional smoothness conditions on the marginal density of $X$. For other popular regression estimators, such as kernel or orthogonal series ones, a similar inequality can also be proved if we assume that the marginal density of $X$ is as smooth as the regression function.

DEFINITION 3.1. For a fixed parameter $\alpha \geq 0$, fixed positive parameters $c_0, r_0, C_0, \beta, L, \mu_{\max} > \mu_{\min} > 0$ and a fixed compact $\mathcal{C} \subset \mathbf{R}^d$, let $\mathcal{P}_\Sigma$ denote the class of all probability distributions $P$ on $\mathcal{Z}$ such that:

(i) the margin Assumption (MA) is satisfied,
(ii) the regression function $\eta$ belongs to the Hölder class $\Sigma(\beta, L, \mathbf{R}^d)$,
(iii) the strong density assumption on $P_X$ is satisfied.

Lemma 3.1 and (3.8) immediately imply the next result.

THEOREM 3.3. *For any $n \geq 1$ the excess risk of the plug-in classifier $\hat{f}_n^* = \mathbb{1}_{\{\hat{\eta}_n^* \geq 1/2\}}$ with bandwidth $h = n^{-1/(2\beta+d)}$ satisfies*

$$\sup_{P \in \mathcal{P}_\Sigma} \{\mathbf{E} R(\hat{f}_n^*) - R(f^*)\} \leq C n^{-\beta(1+\alpha)/(2\beta+d)},$$

*where the constant $C > 0$ depends only on $\alpha$, $C_0$, $C_1$ and $C_2$.*

For $\alpha\beta > d/2$ the convergence rate $n^{-\beta(1+\alpha)/(2\beta+d)}$ obtained in Theorem 3.3 is a *fast rate*, i.e., it is faster than $n^{-1/2}$. Furthermore, it is a *superfast rate* (i.e., is faster than $n^{-1}$) for $\alpha\beta > d$. We must note that if this

FAST LEARNING RATES FOR PLUG-IN CLASSIFIERS 11condition is satisfied, the class $\mathcal{P}_\Sigma$ is rather poor, and thus super-fast rates can occur only for very particular joint distributions of $(X, Y)$. Intuitively this is clear. Indeed, to have a very smooth regression function $\eta$ (i.e., very large $\beta$) implies that when $\eta$ hits the level $1/2$, it cannot "take off" from this level too abruptly. As a consequence, when the density of the distribution $P_X$ is bounded away from 0 in the vicinity of the hitting point, the margin assumption cannot be satisfied for large $\alpha$ since this assumption puts an upper bound on the "time spent" by the regression function near $1/2$. So, $\alpha$ and $\beta$ cannot be simultaneously very large. It can be shown that the case of simultaneously large $\alpha$ and $\beta$ is essentially described by the condition $\alpha\beta > d$.

To be more precise, observe first that $\mathcal{P}_\Sigma$ is not empty for $\alpha\beta > d$, so that super-fast rates can effectively occur. Examples of laws $P \in \mathcal{P}_\Sigma$ under this condition can be easily given, such as the one with $P_X$ equal to the uniform distribution on a ball centered at 0 in $\mathbf{R}^d$, and the regression function defined by $\eta(x) = 1/2 - C\|x\|^2$ with an appropriate $C > 0$. Clearly, $\eta$ belongs to Hölder classes with arbitrarily large $\beta$ and Assumption (MA) is satisfied with $\alpha = d/2$. Thus, for $d \geq 3$ and $\beta$ large enough super-fast rates can occur. Note that in this example the decision set $\{x : \eta(x) \geq 1/2\}$ has Lebesgue measure 0 in $\mathbf{R}^d$. It turns out that such a condition is necessary to achieve classification with super-fast rates when the Hölder classes of regression functions are considered.

To explain this we need a definition. We will say that $\eta$ crosses the level $1/2$ at a point $x_0 \in \mathbf{R}^d$ if for any $r > 0$ there exist $x_-$ and $x_+$ in $\mathcal{B}(x_0, r)$ such that $\eta(x_-) < 1/2$ and $\eta(x_+) > 1/2$.

PROPOSITION 3.4. *If $\alpha(1 \wedge \beta) > 1$ there is no distribution $P \in \mathcal{P}_\Sigma$ such that the regression function $\eta$ associated with $P$ crosses $1/2$ in the interior of the support of $P_X$.*

Proof of this proposition is given in [2].

Note that the condition $\alpha(1 \wedge \beta) > 1$ appearing in Proposition 3.4 is equivalent to $\frac{\beta(1+\alpha)}{2\beta+d} > \frac{(2\beta)\vee(\beta+1)}{2\beta+d}$, which is necessary to have super-fast rates. As a simple consequence, in this context, super-fast rates cannot occur when the regression function crosses $1/2$ in the interior of the support.

The following lower bound shows optimality of the rate of convergence for the Hölder classes obtained in Theorem 3.3.

THEOREM 3.5. *Let $d \geq 1$ be an integer, and let $L, \beta, \alpha$ be positive constants such that $\alpha\beta \leq d$. Then there exists a constant $C > 0$ such that for any $n \geq 1$ and any classifier $\hat{f}_n : \mathcal{Z}^n \to \mathcal{F}$, we have*

$$\sup_{P \in \mathcal{P}_\Sigma} \{\mathbf{E}R(\hat{f}_n) - R(f^*)\} \geq C n^{-\beta(1+\alpha)/(2\beta+d)}.$$

The proof is given in Section 6.2.

Note that the lower bound of Theorem 3.5 does not cover the case of super-fast rates ($\alpha\beta > d$).

Finally, we discuss the case where "$\alpha = \infty$," which means that there exists $t_0 > 0$ such that

$$P_X(0 < |\eta(X) - 1/2| \leq t_0) = 0. \tag{3.9}$$

This is a very favorable situation for classification. The rates of convergence of the ERM type classifiers under (3.9) are, of course, faster than under Assumption (MA) with $\alpha < \infty$ (cf. [15]), but they are not faster than $n^{-1}$. Indeed, Massart and Nédélec [15] provide a lower bound showing that, even if Assumption (CAD) is replaced by the very strong assumption that the true decision set belongs to a VC-class (note that both assumptions are naturally linked to the study of the ERM type classifiers), the best achievable rate is of order $(\log n)/n$. We show now that for the plug-in classifiers much faster rates can be attained. Specifically, if the regression function $\eta$ has some (arbitrarily low) Hölder smoothness $\beta$ and (3.9) holds, the rate of convergence is exponential in $n$. To show this, we first state a simple lemma which is valid for any plug-in classifier $\hat{f}_n$.

LEMMA 3.6. *Let assumption* (3.9) *be satisfied, and let $\hat{\eta}_n$ be an estimator of the regression function $\eta$. Then for the plug-in classifier $\hat{f}_n = \mathbb{1}_{\{\hat{\eta}_n \geq 1/2\}}$ we have*

$$\mathbf{E}R(\hat{f}_n) - R(f^*) \leq \mathbf{P}(|\hat{\eta}_n(X) - \eta(X)| > t_0).$$

PROOF. Following an argument similar to the proof of Lemma 3.1 and using condition (3.9), we get

$$\begin{aligned}
\mathbf{E}R(\hat{f}_n) - R(f^*) &\leq 2t_0 P_X(0 < |\eta(X) - 1/2| \leq t_0) \\
&\quad + \mathbf{E}(|2\eta(X) - 1|\mathbb{1}_{\{\hat{f}_n(X) \neq f^*(X)\}} \mathbb{1}_{\{|\eta(X) - 1/2| > t_0\}}) \\
&= \mathbf{E}(|2\eta(X) - 1|\mathbb{1}_{\{\hat{f}_n(X) \neq f^*(X)\}} \mathbb{1}_{\{|\eta(X) - 1/2| > t_0\}}) \\
&\leq \mathbf{P}(|\hat{\eta}_n(X) - \eta(X)| > t_0). \qquad \square
\end{aligned}$$

Lemma 3.6 and Theorem 3.2 immediately imply that, under assumption (3.9), the rate of convergence of the plug-in classifier $\hat{f}_n^* = \mathbb{1}_{\{\hat{\eta}_n^* \geq 1/2\}}$ with a small enough fixed (independent of $n$) bandwidth $h$ is exponential. To state the result, we denote by $\mathcal{P}_{\Sigma,\infty}$ the class of probability distributions $P$ defined in the same way as $\mathcal{P}_\Sigma$, with the only difference being that in Definition 3.1 the margin Assumption (MA) is replaced by condition (3.9).

PROPOSITION 3.7. *There exists a fixed (independent of $n$) $h > 0$ such that for any $n \geq 1$ the excess risk of the plug-in classifier $\hat{f}_n^* = \mathbb{1}_{\{\hat{\eta}_n^* \geq 1/2\}}$ with bandwidth $h$ satisfies*

$$\sup_{P \in \mathcal{P}_{\Sigma,\infty}} \{\mathbf{E}R(\hat{f}_n^*) - R(f^*)\} \leq C_4 \exp(-C_5 n),$$

*where the constants $C_4, C_5 > 0$ depend only on $t_0$, $\beta$, $d$, $L$, $c_0$, $r_0$, $\mu_{\min}, \mu_{\max}$, and on the kernel $K$.*

PROOF. Use Lemma 3.6, choose $h > 0$ such that $h < \min(r_0/c, (t_0/C_3)^{1/\beta})$ and apply (3.7) with $\delta = t_0$. □

Koltchinskii and Beznosova [12] prove a result on exponential rates for the plug-in classifier with a penalized regression estimator in place of the local polynomial one that we use here. Their result is stated under a less general condition than Proposition 3.7, in the sense that they consider only the Lipschitz class of regression functions $\eta$, while in Proposition 3.7 the Hölder smoothness $\beta$ can be arbitrarily close to 0. Note also that we do not impose any complexity assumption on the decision set. However, the class of distributions $\mathcal{P}_{\Sigma,\infty}$ is quite restricted in a different sense. Indeed, for such distributions condition (3.9) should be compatible with the assumption that $\eta$ belongs to a Hölder class. A sufficient condition for this is the existence of a band or a "corridor" of zero $P_X$-measure separating the sets $\{x : \eta(x) > 1/2\}$ and $\{x : \eta(x) < 1/2\}$. We believe that this condition is close to the necessary one.

**4. Optimal learning rates without the strong density assumption.** In this section we show that if $P_X$ does not admit a density bounded away from zero on its support the rates of classification are slower than those obtained in Section 3. In particular, super-fast rates, that is, the rates faster than $n^{-1}$, cannot be achieved. Introduce the following class of probability distributions.

DEFINITION 4.1. For a fixed parameter $\alpha \geq 0$, fixed positive parameters $c_0, r_0, C_0, \beta, L, \mu_{\max} > 0$ and a fixed compact $\mathcal{C} \subset \mathbf{R}^d$, let $\mathcal{P}'_\Sigma$ denote the class of all probability distributions $P$ on $\mathcal{Z}$ such that:

(i) the margin Assumption (MA) is satisfied,
(ii) the regression function $\eta$ belongs to the Hölder class $\Sigma(\beta, L, \mathbf{R}^d)$,
(iii) the mild density assumption on $P_X$ is satisfied.

In this section we mainly assume that the distribution $P$ of $(X, Y)$ belongs to $\mathcal{P}'_\Sigma$, but we also consider larger classes of distributions satisfying the margin Assumption (MA) and the complexity Assumption (CAR).

Clearly, $\mathcal{P}_\Sigma \subset \mathcal{P}'_\Sigma$. The only difference between $\mathcal{P}'_\Sigma$ and $\mathcal{P}_\Sigma$ is that for $\mathcal{P}'_\Sigma$ the marginal density of $X$ is not bounded away from zero. The optimal rates for $\mathcal{P}'_\Sigma$ are slower than for $\mathcal{P}_\Sigma$. Indeed, we have the following lower bound for the excess risk.

THEOREM 4.1. *Let $d \geq 1$ be an integer, and let $L, \beta, \alpha$ be positive constants. Then there exists a constant $C > 0$ such that for any $n \geq 1$ and any classifier $\hat{f}_n : \mathcal{Z}^n \to \mathcal{F}$ we have*

$$\sup_{P \in \mathcal{P}'_\Sigma} \{\mathbf{E}R(\hat{f}_n) - R(f^*)\} \geq C n^{-(1+\alpha)\beta/((2+\alpha)\beta+d)}.$$

The proof is given in Section 6.2.

In particular, when $\alpha = d/\beta$, we get the slow convergence rate $1/\sqrt{n}$, instead of the fast rate $n^{-(\beta+d)/(2\beta+d)}$ obtained in Theorem 3.3 under the strong density assumption. Nevertheless, the lower bound can still approach $n^{-1}$ as the margin parameter $\alpha$ tends to $\infty$.

We now show that the rate of convergence given in Theorem 4.1 is optimal in the sense that there exist estimators that achieve this rate. This will be obtained as a consequence of a general upper bound for the excess risk of classifiers over a larger set $\mathcal{P}$ of distributions than $\mathcal{P}'_\Sigma$.

Fix a Lebesgue measurable set $\mathcal{C} \subseteq \mathbf{R}^d$ and a value $1 \leq p \leq \infty$. Let $\Sigma$ be a class of regression functions $\eta$ on $\mathbf{R}^d$ such that Assumption (CAR) is satisfied where the $\varepsilon$-entropy is taken w.r.t. the $L_p(\mathcal{C}, P_X)$ norm. Then for every $\varepsilon > 0$ there exists an $\varepsilon$-net $\mathcal{N}_\varepsilon$ on $\Sigma$ w.r.t. this norm such that

$$\log(\operatorname{card} \mathcal{N}_\varepsilon) \leq A' \varepsilon^{-\rho},$$

where $A'$ is a constant. Consider the empirical risk

$$R_n(f) = \frac{1}{n} \sum_{i=1}^n \mathbb{1}_{\{f(X_i) \neq Y_i\}}, \qquad f \in \mathcal{F},$$

and set

$$\varepsilon_n = \varepsilon_n(\alpha, \rho, p) \triangleq \begin{cases} n^{-1/(2+\alpha+\rho)}, & \text{if } p = \infty, \\ n^{-(p+\alpha)/((2+\alpha)p+\rho(p+\alpha))}, & \text{if } 1 \leq p < \infty. \end{cases}$$

Define a sieve estimator $\hat{\eta}_n^S$ of the regression function $\eta$ by the relation

(4.1) $$\hat{\eta}_n^S \in \operatorname*{Arg\,min}_{\bar{\eta} \in \mathcal{N}_{\varepsilon_n}} R_n(f_{\bar{\eta}}),$$

where $f_{\bar{\eta}}(x) = \mathbb{1}_{\{\bar{\eta}(x) \geq 1/2\}}$, and consider the classifier $\hat{f}_n^S = \mathbb{1}_{\{\hat{\eta}_n^S \geq 1/2\}}$. Note that $\hat{f}_n^S$ can be viewed as a "hybrid" plug-in/ ERM procedure: the ERM is performed on a set of plug-in rules corresponding to a grid on the class of regression functions $\eta$.

THEOREM 4.2. *Let $\mathcal{P}$ be a set of probability distributions $P$ on $\mathcal{Z}$ such that:*

(i) *the margin Assumption* (MA) *is satisfied,*

(ii) *the regression function $\eta$ belongs to a class $\Sigma$ which satisfies the complexity Assumption* (CAR) *with the $\varepsilon$-entropy taken w.r.t. the $L_p(\mathcal{C}, P_X)$ norm for some $1 \leq p \leq \infty$.*

*Consider the classifier $\hat{f}_n^S = \mathbb{1}_{\{\hat{\eta}_n^S \geq 1/2\}}$. Then for any $n \geq 1$ we have*

$$\sup_{P \in \mathcal{P}} \{\mathbf{E}R(\hat{f}_n^S) - R(f^*)\}$$

(4.2)
$$\leq \begin{cases} Cn^{-(1+\alpha)/(2+\alpha+\rho)}, & \text{if } p = \infty, \\ Cn^{-(1+\alpha)p/((2+\alpha)p+\rho(p+\alpha))}, & \text{if } 1 \leq p < \infty. \end{cases}$$

The proof is given in Section 6.3.

Theorem 4.2 allows one to get fast classification rates without any density assumption on $P_X$. Namely, define the following class of distributions $P$ of $(X, Y)$.

DEFINITION 4.2. For fixed parameters $\alpha \geq 0$, $C_0 > 0, \beta > 0, L > 0$, and for a fixed compact $\mathcal{C} \subset \mathbf{R}^d$, let $\mathcal{P}_\Sigma^0$ denote the class of all probability distributions $P$ on $\mathcal{Z}$ such that:

(i) the margin Assumption (MA) is satisfied,
(ii) the regression function $\eta$ belongs to the Hölder class $\Sigma(\beta, L, \mathbf{R}^d)$,
(iii) for all $P \in \mathcal{P}$ the supports of marginal distributions $P_X$ are included in $\mathcal{C}$.

If $\mathcal{C}$ is a compact, the estimates of $\varepsilon$-entropies of Hölder classes $\Sigma(\beta, L, \mathbf{R}^d)$ in the $L_\infty(\mathcal{C}, \lambda)$ norm where $\lambda$ is the Lebesgue measure on $\mathbf{R}^d$ are obtained by Kolmogorov and Tihomirov [10], and they yield Assumption (CAR) with $\rho = d/\beta$. Therefore, from (4.2) with $p = \infty$ we easily get the following upper bound.

THEOREM 4.3. *Let $d \geq 1$ be an integer, and let $L, \beta$ and $\alpha$ be positive constants. For any $n \geq 1$ the classifier $\hat{f}_n^S = \mathbb{1}_{\{\hat{\eta}_n^S \geq 1/2\}}$ defined by* (4.1) *with $p = \infty$ satisfies*

$$\sup_{P \in \mathcal{P}_\Sigma^0} \{\mathbf{E}R(\hat{f}_n^S) - R(f^*)\} \leq Cn^{-(1+\alpha)\beta/((2+\alpha)\beta+d)}$$

*with some constant $C > 0$ depending only on $\alpha$, $\beta$, $d$, $L$ and $C_0$.*

Since $\mathcal{P}'_\Sigma \subset \mathcal{P}_\Sigma^0$, Theorems 3.5 and 4.3 show that $n^{-(1+\alpha)\beta/((2+\alpha)\beta+d)}$ is the optimal rate of convergence of the excess risk on the class of distributions $\mathcal{P}_\Sigma^0$.

**5. Comparison lemmas.** In this section we give some useful inequalities between the risks of plug-in classifiers and the $L_p$ risks of the corresponding regression estimators under the margin Assumption (MA). These inequalities will be helpful in the proofs. They also illustrate a connection between the two complexity Assumptions (CAR) and (CAD) defined in the Introduction and allow one to compare our study of plug-in estimators with that given by Yang [26], who considered the case $\alpha = 0$ (no margin assumption), as well as with the developments in [3] and [6].

Throughout this section $\bar{\eta}$ is a Borel function on $\mathbf{R}^d$ and

$$\bar{f}(x) = \mathbb{1}_{\{\bar{\eta}(x) \geq 1/2\}}.$$

For $1 \leq p \leq \infty$ we denote by $\|\cdot\|_p$ the $L_p(\mathbf{R}^d, P_X)$ norm. We first state some comparison inequalities for the $L_\infty$ norm.

LEMMA 5.1. *For any distribution $P$ of $(X,Y)$ satisfying Assumption* (MA) *we have*

(5.1) $$R(\bar{f}) - R(f^*) \leq 2C_0 \|\bar{\eta} - \eta\|_\infty^{1+\alpha}$$

*and*

(5.2) $$P_X(\bar{f}(X) \neq f^*(X), \eta(X) \neq 1/2) \leq C_0 \|\bar{\eta} - \eta\|_\infty^\alpha.$$

PROOF. To show (5.1) note that

$$\begin{aligned}
R(\bar{f}) - R(f^*) &= \mathbf{E}(|2\eta(X) - 1|\mathbb{1}_{\{\bar{f}(X) \neq f^*(X)\}}) \\
&\leq 2\mathbf{E}(|\eta(X) - \tfrac{1}{2}|\mathbb{1}_{\{0 < |\eta(X) - 1/2| \leq |\eta(X) - \bar{\eta}(X)|\}}) \\
&\leq 2\|\eta - \bar{\eta}\|_\infty P_X(0 < |\eta(X) - \tfrac{1}{2}| \leq \|\eta - \bar{\eta}\|_\infty) \\
&\leq 2C_0 \|\eta - \bar{\eta}\|_\infty^{1+\alpha}.
\end{aligned}$$

Similarly,

$$\begin{aligned}
P_X(\bar{f}(X) \neq f^*(X), \eta(X) \neq 1/2) &\leq P_X(0 < |\eta(X) - \tfrac{1}{2}| \leq |\eta(X) - \bar{\eta}(X)|) \\
&\leq P_X(0 < |\eta(X) - \tfrac{1}{2}| \leq \|\eta - \bar{\eta}\|_\infty) \\
&\leq C_0 \|\eta - \bar{\eta}\|_\infty^\alpha. \qquad \square
\end{aligned}$$

REMARK 5.1. Lemma 5.1 offers an easy way to obtain the result of Theorem 3.3 in a slightly less precise form, with an extra logarithmic factor in the rate. In fact, under the strong density assumption, there exist nonparametric estimators $\hat{\eta}_n$ (e.g., suitably chosen local polynomial estimators) such that

$$\mathbf{E}(\|\hat{\eta}_n - \eta\|_\infty^q) \leq C\left(\frac{\log n}{n}\right)^{q\beta/(2\beta+d)} \qquad \forall q > 0,$$

uniformly over $\eta \in \Sigma(\beta, L, \mathbf{R}^d)$ (see, e.g., [18]). Taking here $q = 1 + \alpha$ and applying the comparison inequality (5.1) we immediately get that the plug-in classifier $\hat{f}_n = \mathbb{1}_{\{\hat{\eta}_n \geq 1/2\}}$ has excess risk $\mathcal{E}(\hat{f}_n)$ of the order $(n/\log n)^{-\beta(1+\alpha)/(2\beta+d)}$.

Another immediate application of Lemma 5.1 is to get lower bounds on the risks of regression estimators in the $L_\infty$ norm from the corresponding lower bounds on the excess risks of classifiers (cf. Theorems 3.5 and 4.1). But here again we lose a logarithmic factor required for the best bounds.

Inequality (5.2) serves to compare the measure of symmetric difference distance between the decision sets with the $L_\infty$ distance between the corresponding regression functions. In fact, if $P_X(\eta(X) = 1/2) = 0$, inequality (5.2) reads as $d_\triangle(G, G^*) \leq C_0 \|\bar{\eta} - \eta\|_\infty^\alpha$, where $G = \{x : \bar{f}(x) = 1\}$. Thus, the $d_\triangle$-convergence rates for estimation of the true decision set $G^*$ can be obtained from the $L_\infty$ rates of the corresponding regression estimators.

We now consider the comparison inequalities for $L_p$ norms with $1 \leq p < \infty$.

LEMMA 5.2. *For any $1 \leq p < \infty$ and any distribution $P$ of $(X,Y)$ satisfying Assumption* (MA) *with $\alpha > 0$ we have*

$$(5.3) \qquad R(\bar{f}) - R(f^*) \leq C_1(\alpha, p) \|\bar{\eta} - \eta\|_p^{p(1+\alpha)/(p+\alpha)},$$

*where $C_1(\alpha, p) = 2(\alpha + p)p^{-1}(\frac{p}{\alpha})^{\alpha/(\alpha+p)} C_0^{(p-1)/(\alpha+p)}$. In particular,*

$$(5.4) \quad R(\bar{f}) - R(f^*) \leq C_1(\alpha, 2) \left( \int [\bar{\eta}(x) - \eta(x)]^2 P_X(dx) \right)^{(1+\alpha)/(2+\alpha)}.$$

PROOF. For any $t > 0$ we have

$$\begin{aligned}
R(\bar{f}) - R(f^*) &= \mathbf{E}[|2\eta(X) - 1|\mathbb{1}_{\{\bar{f}(X) \neq f^*(X)\}}] \\
&= 2\mathbf{E}[|\eta(X) - 1/2|\mathbb{1}_{\{\bar{f}(X) \neq f^*(X)\}}\mathbb{1}_{\{0 < |\eta(X) - 1/2| \leq t\}}] \\
&\quad + 2\mathbf{E}[|\eta(X) - 1/2|\mathbb{1}_{\{\bar{f}(X) \neq f^*(X)\}}\mathbb{1}_{\{|\eta(X) - 1/2| > t\}}] \\
(5.5) \quad &\leq 2\mathbf{E}[|\eta(X) - \bar{\eta}(X)|\mathbb{1}_{\{0 < |\eta(X) - 1/2| \leq t\}}] \\
&\quad + 2\mathbf{E}[|\eta(X) - \bar{\eta}(X)|\mathbb{1}_{\{|\eta(X) - \bar{\eta}(X)| > t\}}] \\
&\leq 2\|\eta - \bar{\eta}\|_p [P_X(0 < |\eta(X) - 1/2| \leq t)]^{(p-1)/p} \\
&\quad + \frac{2\|\eta - \bar{\eta}\|_p^p}{t^{p-1}},
\end{aligned}$$

by the Hölder and Markov inequalities. So, for any $t > 0$, introducing $E \stackrel{\triangle}{=} \|\eta - \bar{\eta}\|_p$ and using Assumption (MA) to bound the probability in (5.5), we

obtain

$$R(\bar{f}) - R(f^*) \leq 2\bigg(C_0^{(p-1)/p} t^{\alpha(p-1)/p} E + \frac{E^p}{t^{p-1}}\bigg).$$

Minimizing over $t$ the RHS of this inequality we get (5.3). □

If the regression function $\eta$ belongs to the Hölder class $\Sigma(\beta, L, \mathbf{R}^d)$ there exist estimators $\hat{\eta}_n$ such that, uniformly over the class,

(5.6) $\qquad \mathbf{E}\{[\hat{\eta}_n(X) - \eta(X)]^2\} \leq Cn^{-2\beta/(2\beta+d)}$

for some constant $C > 0$. This has been shown by Stone [18] under the additional strong density assumption and by Yang [26] with no assumption on $P_X$. Using (5.6) and (5.4) we get that the excess risk of the corresponding plug-in classifier $\hat{f}_n = \mathbb{1}_{\{\hat{\eta}_n \geq 1/2\}}$ admits a bound of the order $n^{-(2\beta/(2\beta+d))((1+\alpha)/(2+\alpha))}$, which is suboptimal when $\alpha \neq 0$ (cf. Theorems 4.2, 4.3). In other words, under the margin assumption, Lemma 5.2 is not the right tool to analyze the convergence rate of plug-in classifiers. On the contrary, when no margin assumption is imposed (i.e., $\alpha = 0$ in our notation), inequality (1.3), which is a version of (5.4) for $\alpha = 0$, is precise enough to give the optimal rate of classification [26].

Another way to obtain (5.4) is to use Bartlett, Jordan and McAuliffe [3]. It is enough to apply their Theorem 10 with (in their notation) $\phi(t) = (1-t)^2$, $\psi(t) = t^2$, and to note that for this choice of $\phi$ we have $R_\phi(\bar{\eta}) - R_\phi^* = \|\eta - \bar{\eta}\|_2^2$. Blanchard, Lugosi and Vayatis [6] used the result of Bartlett, Jordan and McAuliffe [3] to prove fast rates of the order $n^{-2(1+\alpha)/(3(2+\alpha))}$ for a boosting procedure over the class of regression functions $\eta$ of bounded variation in dimension $d = 1$. Note that the same rates can be obtained for other plug-in classifiers using (5.4). Indeed, if $\eta$ is of bounded variation, there exist estimators of $\eta$ converging with the mean squared $L_2$ rate $n^{-2/3}$ (cf. [9, 16, 23, 27]), and thus application of (5.4) immediately yields the rate $n^{-2(1+\alpha)/(3(2+\alpha))}$ for the corresponding plug-in rule. However, Theorem 4.2 shows that this is not an optimal rate [here again we observe that inequality (5.4) fails to establish the optimal properties of plug-in classifiers]. In fact, let $d = 1$ and let the assumptions of Theorem 4.2 be satisfied, where instead of assumption (ii) we use a particular case: $\eta$ belongs to a class of functions on $[0, 1]$ whose total variation is bounded by a constant $L < \infty$. It follows from [4] that Assumption (CAR) for this class is satisfied with $\rho = 1$ for any $1 \leq p < \infty$. Hence, we can apply (4.2) of Theorem 4.2 to find that

(5.7) $\qquad \sup_{P \in \mathcal{P}} \{\mathbf{E}R(\hat{f}_n^S) - R(f^*)\} \leq Cn^{-(1+\alpha)p/((2+\alpha)p+(p+\alpha))}$

for the corresponding class $\mathcal{P}$. If $p > 2$ (recall that the value $p \in [1, \infty)$ is chosen by the statistician), the rate in (5.7) is faster than $n^{-2(1+\alpha)/(3(2+\alpha))}$ obtained under the same conditions by Blanchard, Lugosi and Vayatis [6].

**6. Proofs.**

6.1. *Proof of Theorem* 3.2. Consider a distribution $P$ in $\mathcal{P}_\Sigma$. Let $A$ be the support of $P_X$. Fix $x \in A$ and $\delta > 0$. Consider the matrix $\bar{B} \triangleq (\bar{B}_{s_1,s_2})_{|s_1|,|s_2| \leq \lfloor \beta \rfloor}$ with elements $\bar{B}_{s_1,s_2} \triangleq \int_{\mathbf{R}^d} u^{s_1+s_2} K(u) \mu(x+hu)\,du$. The smallest eigenvalue $\lambda_{\bar{B}}$ of $\bar{B}$ satisfies

$$\lambda_{\bar{B}} = \min_{\|W\|=1} W^T \bar{B} W$$

(6.1)
$$\geq \min_{\|W\|=1} W^T BW + \min_{\|W\|=1} W^T (\bar{B} - B) W$$

$$\geq \min_{\|W\|=1} W^T BW - \sum_{|s_1|,|s_2| \leq \lfloor \beta \rfloor} |\bar{B}_{s_1,s_2} - B_{s_1,s_2}|.$$

Let $A_n \triangleq \{u \in \mathbf{R}^d : \|u\| \leq c; x + hu \in A\}$ where $c$ is the constant appearing in (3.3). Using (3.3), for any vector $W$ satisfying $\|W\| = 1$, we obtain

$$W^T BW = \int_{\mathbf{R}^d} \left( \sum_{|s| \leq \lfloor \beta \rfloor} W_s u^s \right)^2 K(u) \mu(x+hu)\,du$$

$$\geq c\mu_{\min} \int_{A_n} \left( \sum_{|s| \leq \lfloor \beta \rfloor} W_s u^s \right)^2 du.$$

By assumption of the theorem, $ch \leq r_0$. Since the support of the marginal distribution is $(c_0, r_0)$-regular, we get

$$\lambda[A_n] \geq h^{-d} \lambda[\mathcal{B}(x, ch) \cap A] \geq c_0 h^{-d} \lambda[\mathcal{B}(x, ch)] \geq c_0 v_d c^d,$$

where $v_d \triangleq \lambda[\mathcal{B}(0,1)]$ is the volume of the unit ball and $c_0 > 0$ is the constant introduced in the definition (2.1) of the $(c_0, r_0)$-regular set.

Let $\mathcal{A}$ denote the class of all compact subsets of $\mathcal{B}(0, c)$ having Lebesgue measure $c_0 v_d c^d$. Using the previous displays we obtain

(6.2) $$\min_{\|W\|=1} W^T BW \geq c\mu_{\min} \min_{\|W\|=1; S \in \mathcal{A}} \int_S \left( \sum_{|s| \leq \lfloor \beta \rfloor} W_s u^s \right)^2 du \triangleq 2\mu_0.$$

By the compactness argument, the minimum in (6.2) exists and is strictly positive.

For $i = 1, \ldots, n$ and any multi-indices $s_1, s_2$ such that $|s_1|, |s_2| \leq \lfloor \beta \rfloor$, define

$$T_i \triangleq \frac{1}{h^d} \left( \frac{X_i - x}{h} \right)^{s_1+s_2} K\left( \frac{X_i - x}{h} \right) - \int_{\mathbf{R}^d} u^{s_1+s_2} K(u) \mu(x+hu)\,du.$$

We have $\mathbf{E}T_i = 0$, $|T_i| \leq h^{-d}\sup_{u\in\mathbf{R}^d}(1+\|u\|^{2\beta})K(u) \triangleq \kappa_1 h^{-d}$ and the following bound for the variance of $T_i$:

$$\mathbf{Var}\,T_i \leq \frac{1}{h^{2d}}\mathbf{E}\left(\frac{X_i-x}{h}\right)^{2s_1+2s_2}K^2\left(\frac{X_i-x}{h}\right)$$

$$= \frac{1}{h^d}\int_{\mathbf{R}^d} u^{2s_1+2s_2}K^2(u)\mu(x+hu)\,du$$

$$\leq \frac{\mu_{\max}}{h^d}\int_{\mathbf{R}^d}(1+\|u\|^{4\beta})K^2(u)\,du \triangleq \frac{\kappa_2}{h^d}.$$

Using Bernstein's inequality, for any $\varepsilon > 0$, we have

$$P^{\otimes n}(|\bar{B}_{s_1,s_2} - B_{s_1,s_2}| > \varepsilon) = P^{\otimes n}\left(\left|\frac{1}{n}\sum_{i=1}^n T_i\right| > \varepsilon\right) \leq 2\exp\left\{-\frac{nh^d\varepsilon^2}{2\kappa_2 + 2\kappa_1\varepsilon/3}\right\}.$$

This and (6.1) and (6.2) imply that

$$(6.3)\qquad P^{\otimes n}(\lambda_{\bar{B}} \leq \mu_0) \leq 2M^2 \exp(-Cnh^d),$$

where $M^2$ is the number of elements of the matrix $\bar{B}$. Assume in what follows that $n$ is large enough so that $\mu_0 > (\log n)^{-1}$. Then for $\lambda_{\bar{B}} > \mu_0$ we have $|\hat{\eta}_n^*(x) - \eta(x)| \leq |\hat{\eta}_n^{\mathrm{LP}}(x) - \eta(x)|$. Therefore,

$$(6.4)\quad\begin{aligned}P^{\otimes n}(|\hat{\eta}_n^*(x) - \eta(x)| \geq \delta) &\leq P^{\otimes n}(\lambda_{\bar{B}} \leq \mu_0)\\ &\quad + P^{\otimes n}(|\hat{\eta}_n^{\mathrm{LP}}(x) - \eta(x)| \geq \delta, \lambda_{\bar{B}} > \mu_0).\end{aligned}$$

We now evaluate the second probability on the right-hand side of (6.4). For $\lambda_{\bar{B}} > \mu_0$ we have $\hat{\eta}_n^{\mathrm{LP}}(x) = U^T(0)Q^{-1}V$ [where $V$ is given by (2.3)]. Introduce the matrix $Z \triangleq (Z_{i,s})_{1\leq i\leq n, |s|\leq \lfloor\beta\rfloor}$ with elements

$$Z_{i,s} \triangleq (X_i - x)^s\sqrt{K\left(\frac{X_i-x}{h}\right)}.$$

The $s$th column of $Z$ is denoted by $Z_s$, and we introduce

$$Z^{(\eta)} \triangleq \sum_{|s|\leq \lfloor\beta\rfloor} \frac{\eta^{(s)}(x)}{s!}Z_s.$$

Since $Q = Z^T Z$, we get

$$\forall |s| \leq \lfloor\beta\rfloor: \qquad U^T(0)Q^{-1}Z^T Z_s = \mathbb{1}_{\{s=(0,\ldots,0)\}},$$

hence $U^T(0)Q^{-1}Z^T Z^{(\eta)} = \eta(x)$. So we can write

$$\hat{\eta}_n^{\mathrm{LP}}(x) - \eta(x) = U^T(0)Q^{-1}(V - Z^T Z^{(\eta)}) = U^T(0)\bar{B}^{-1}\mathbf{a},$$

where $\mathbf{a} \triangleq \frac{1}{nh^d} H(V - Z^T Z^{(\eta)}) \in \mathbf{R}^M$ and $H$ is the diagonal matrix $H \triangleq (H_{s_1,s_2})_{|s_1|,|s_2| \leq \lfloor \beta \rfloor}$ with $H_{s_1,s_2} \triangleq h^{-s_1} \mathbb{1}_{\{s_1 = s_2\}}$. For $\lambda_{\bar{B}} > \mu_0$ we get

$$(6.5) \quad |\hat{\eta}_n^{\mathrm{LP}}(x) - \eta(x)| \leq \|\bar{B}^{-1} \mathbf{a}\| \leq \lambda_{\bar{B}}^{-1} \|\mathbf{a}\| \leq \mu_0^{-1} \|\mathbf{a}\| \leq \mu_0^{-1} M \max_s |a_s|,$$

where $a_s$ are the components of the vector $\mathbf{a}$ given by

$$a_s = \frac{1}{nh^d} \sum_{i=1}^n [Y_i - \eta_x(X_i)] \left(\frac{X_i - x}{h}\right)^s K\left(\frac{X_i - x}{h}\right).$$

Define

$$T_i^{(s,1)} \triangleq \frac{1}{h^d} [Y_i - \eta(X_i)] \left(\frac{X_i - x}{h}\right)^s K\left(\frac{X_i - x}{h}\right),$$

$$T_i^{(s,2)} \triangleq \frac{1}{h^d} [\eta(X_i) - \eta_x(X_i)] \left(\frac{X_i - x}{h}\right)^s K\left(\frac{X_i - x}{h}\right).$$

We have

$$(6.6) \quad |a_s| \leq \left|\frac{1}{n} \sum_{i=1}^n T_i^{(s,1)}\right| + \left|\frac{1}{n} \sum_{i=1}^n [T_i^{(s,2)} - \mathbf{E} T_i^{(s,2)}]\right| + |\mathbf{E} T_i^{(s,2)}|.$$

Note that $\mathbf{E} T_i^{(s,1)} = 0$, $|T_i^{(s,1)}| \leq \kappa_1 h^{-d}$ and

$$\mathbf{Var}\, T_i^{(s,1)} \leq 4^{-1} h^{-d} \int u^{2s} K^2(u) \mu(x + hu)\, du \leq (\kappa_2/4) h^{-d},$$

$$|T_i^{(s,2)} - \mathbf{E} T_i^{(s,2)}| \leq L\kappa_1 h^{\beta-d} + L\kappa_2 h^\beta \leq C h^{\beta-d},$$

$$\mathbf{Var}\, T_i^{(s,2)} \leq h^{-d} L^2 \int h^{2\beta} \|u\|^{2s+2\beta} K^2(u) \mu(x + hu)\, du \leq L^2 \kappa_2 h^{2\beta-d}.$$

Using Bernstein's inequality, for any $\varepsilon_1, \varepsilon_2 > 0$, we obtain

$$P^{\otimes n}\left(\left|\frac{1}{n} \sum_{i=1}^n T_i^{(s,1)}\right| \geq \varepsilon_1\right) \leq 2 \exp\left\{-\frac{nh^d \varepsilon_1^2}{\kappa_2/2 + 2\kappa_1 \varepsilon_1/3}\right\}$$

and

$$P^{\otimes n}\left(\left|\frac{1}{n} \sum_{i=1}^n [T_i^{(s,2)} - \mathbf{E} T_i^{(s,2)}]\right| \geq \varepsilon_2\right) \leq 2 \exp\left\{-\frac{nh^d \varepsilon_2^2}{2L^2 \kappa_2 h^{2\beta} + 2C h^\beta \varepsilon_2/3}\right\}.$$

Since also

$$|\mathbf{E} T_i^{(s,2)}| \leq L h^\beta \int \|u\|^{s+\beta} K^2(u) \mu(x + hu)\, du \leq L \kappa_2 h^\beta,$$

we get, using (6.6), that if $3\mu_0^{-1}ML\kappa_2 h^\beta \leq \delta \leq 1$ the following inequality holds:

$$P^{\otimes n}\left(|a_s| \geq \frac{\mu_0\delta}{M}\right) \leq P^{\otimes n}\left(\left|\frac{1}{n}\sum_{i=1}^n T_i^{(s,1)}\right| > \frac{\mu_0\delta}{3M}\right)$$
$$+ P^{\otimes n}\left(\left|\frac{1}{n}\sum_{i=1}^n [T_i^{(s,2)} - \mathbf{E}T_i^{(s,2)}]\right| > \frac{\mu_0\delta}{3M}\right)$$
$$\leq 4\exp(-Cnh^d\delta^2).$$

Combining this inequality with (6.3), (6.4) and (6.5), we obtain

(6.7) $\qquad P^{\otimes n}(|\hat\eta_n^*(x) - \eta(x)| \geq \delta) \leq C_1 \exp(-C_2 nh^d \delta^2)$

for $3m^{-1}ML\kappa_2 h^\beta \leq \delta$ (for $\delta > 1$ inequality (6.7) is obvious since $\hat\eta_n^*, \eta$ take values in $[0,1]$). The constants $C_1, C_2$ in (6.7) do not depend on the distribution $P_X$, on its support $A$ and on the point $x \in A$, so that we get (3.7). Now, (3.7) implies (3.8) for $Cn^{-\beta/(2\beta+d)} \leq \delta$, and thus for all $\delta > 0$ (with possibly modified constants $C_1$ and $C_2$).

6.2. *Proof of Theorems* 3.5 *and* 4.1. The proof of both theorems is based on Assouad's lemma (see, e.g., [13], Chapter 2 or [21], Chapter 2). We apply it in a form adapted for the classification problem (Lemma 5.1 in [1]).

For an integer $q \geq 1$ we consider the regular grid on $\mathbf{R}^d$ defined as

$$G_q \triangleq \left\{\left(\frac{2k_1+1}{2q}, \ldots, \frac{2k_d+1}{2q}\right) : k_i \in \{0, \ldots, q-1\}, i = 1, \ldots, d\right\}.$$

Let $n_q(x) \in G_q$ be the closest point to $x \in \mathbf{R}^d$ among points in $G_q$ [we assume uniqueness of $n_q(x)$: if there exist several points in $G_q$ closest to $x$ we define $n_q(x)$ as the one which is closest to 0]. Consider the partition $\mathcal{X}'_1, \ldots, \mathcal{X}'_{q^d}$ of $[0,1]^d$ canonically defined using the grid $G_q$ [$x$ and $y$ belong to the same subset if and only if $n_q(x) = n_q(y)$]. Fix an integer $m \leq q^d$. For any $i \in \{1, \ldots, m\}$, we define $\mathcal{X}_i \triangleq \mathcal{X}'_i$ and $\mathcal{X}_0 \triangleq \mathbf{R}^d \setminus \bigcup_{i=1}^m \mathcal{X}_i$, so that $\mathcal{X}_0, \ldots, \mathcal{X}_m$ form a partition of $\mathbf{R}^d$.

Let $u: \mathbf{R}_+ \to \mathbf{R}_+$ be a nonincreasing infinitely differentiable function such that $u = 1$ on $[0, 1/4]$ and $u = 0$ on $[1/2, \infty)$. One can take, for example, $u(x) = (\int_{1/4}^{1/2} u_1(t)\, dt)^{-1} \int_x^\infty u_1(t)\, dt$ where the infinitely differentiable function $u_1$ is defined as

$$u_1(x) = \begin{cases} \exp\left\{-\dfrac{1}{(1/2-x)(x-1/4)}\right\}, & \text{for } x \in (1/4, 1/2), \\ 0, & \text{otherwise.} \end{cases}$$

Let $\phi: \mathbf{R}^d \to \mathbf{R}_+$ be the function defined as

$$\phi(x) \triangleq C_\phi u(\|x\|),$$

where the positive constant $C_\phi$ is taken small enough to ensure that $|\phi(x') - \phi_x(x')| \leq L\|x' - x\|^\beta$ for any $x, x' \in \mathbf{R}^d$. Thus, $\phi \in \Sigma(\beta, L, \mathbf{R}^d)$.

Define the hypercube $\mathcal{H} = \{\mathbf{P}_{\vec{\sigma}} : \vec{\sigma} = (\sigma_1, \ldots, \sigma_m) \in \{-1, 1\}^m\}$ of probability distributions $\mathbf{P}_{\vec{\sigma}}$ of $(X, Y)$ on $\mathcal{Z} = \mathbf{R}^d \times \{0, 1\}$ as follows.

For any $\mathbf{P}_{\vec{\sigma}} \in \mathcal{H}$ the marginal distribution of $X$ does not depend on $\vec{\sigma}$, and has a density $\mu$ w.r.t. the Lebesgue measure on $\mathbf{R}^d$ defined in the following way. Fix $0 < w \leq m^{-1}$ and a set $A_0$ of positive Lebesgue measure included in $\mathcal{X}_0$ (the particular choice of $A_0$ will be indicated later), and take: (i) $\mu(x) = w/\lambda[\mathcal{B}(0, (4q)^{-1})]$ if $x$ belongs to the ball $\mathcal{B}(z, (4q)^{-1})$ for some $z \in G_q$, (ii) $\mu(x) = (1 - mw)/\lambda[A_0]$ for $x \in A_0$, (iii) $\mu(x) = 0$ for all other $x$.

Next, the distribution of $Y$ given $X$ for $\mathbf{P}_{\vec{\sigma}} \in \mathcal{H}$ is determined by the regression function $\eta_{\vec{\sigma}}(x) = P(Y = 1|X = x)$ that we define as $\eta_{\vec{\sigma}}(x) = \frac{1 + \sigma_j \varphi(x)}{2}$ for any $x \in \mathcal{X}_j$, $j = 1, \ldots, m$, and $\eta_{\vec{\sigma}} \equiv 1/2$ on $\mathcal{X}_0$, where $\varphi(x) \stackrel{\triangle}{=} q^{-\beta}\phi(q[x - n_q(x)])$. We will assume that $C_\phi \leq 1$ to ensure that $\varphi$ and $\eta_{\vec{\sigma}}$ take values in $[0, 1]$.

For any $s \in \mathbf{N}^d$ such that $|s| \leq \lfloor\beta\rfloor$, the partial derivative $D^s\varphi$ exists and $D^s\varphi(x) = q^{|s|-\beta}D^s\phi(q[x - n_q(x)])$. Therefore, for any $i \in \{1, \ldots, m\}$ and any $x, x' \in \mathcal{X}_i$, we have

$$|\varphi(x') - \varphi_x(x')| \leq L\|x - x'\|^\beta.$$

This implies that for any $\vec{\sigma} \in \{-1, 1\}^m$ the function $\eta_{\vec{\sigma}}$ belongs to the Hölder class $\Sigma(\beta, L, \mathbf{R}^d)$.

We now check the margin assumption. Set $x_0 = (\frac{1}{2q}, \ldots, \frac{1}{2q})$. For any $\vec{\sigma} \in \{-1, 1\}^m$ we have

$$\begin{aligned}
\mathbf{P}_{\vec{\sigma}}(0 &< |\eta_{\vec{\sigma}}(X) - 1/2| \leq t) \\
&= m\mathbf{P}_{\vec{\sigma}}(0 < \phi[q(X - x_0)] \leq 2tq^\beta) \\
&= m \int_{\mathcal{B}(x_0, (4q)^{-1})} \mathbb{1}_{\{0 < \phi[q(x-x_0)] \leq 2tq^\beta\}} \frac{w}{\lambda[\mathcal{B}(0, (4q)^{-1})]} \, dx \\
&= \frac{mw}{\lambda[\mathcal{B}(0, 1/4)]} \int_{\mathcal{B}(0, 1/4)} \mathbb{1}_{\{\phi(x) \leq 2tq^\beta\}} \, dx \\
&= mw\mathbb{1}_{\{t \geq C_\phi/(2q^\beta)\}}.
\end{aligned}$$

Therefore, the margin Assumption (MA) is satisfied if $mw = O(q^{-\alpha\beta})$.

According to Lemma 5.1 in [1], for any classifier $\hat{f}_n$ we have

(6.8) $$\sup_{P \in \mathcal{H}} \{\mathbf{E}R(\hat{f}_n) - R(f^*)\} \geq mwb'(1 - b\sqrt{nw})/2,$$

where

$$b \stackrel{\triangle}{=} \left[1 - \left(\int_{\mathcal{X}_1} \sqrt{1 - \varphi^2(x)}\mu_1(x) \, dx\right)^2\right]^{1/2} = C_\phi q^{-\beta},$$

$$b' \stackrel{\triangle}{=} \int_{\mathcal{X}_1} \varphi(x)\mu_1(x)\, dx = C_\phi q^{-\beta},$$

with $\mu_1(x) = \mu(x)/\int_{\mathcal{X}_1} \mu(z)\, dz$.

We now prove Theorem 3.5. Take $q = \lfloor \bar{C} n^{1/(2\beta+d)} \rfloor$, $w = C' q^{-d}$ and $m = \lfloor C'' q^{d-\alpha\beta} \rfloor$ with some positive constants $\bar{C}, C'$ and $C''$ to be chosen, and set $A_0 = [0,1]^d \setminus \bigcup_{i=1}^m \mathcal{X}_i$. The condition $\alpha\beta \leq d$ ensures that the above choice of $m$ is not degenerate: we have $m \geq 1$ for $C''$ large enough. We now prove that $\mathcal{H} \subset \mathcal{P}_\Sigma$ under the appropriate choice of $\bar{C}, C'$ and $C''$. In fact, select these constants so that the triplet $(q, w, m)$ meets the conditions $m \leq q^d$, $0 < w \leq m^{-1}$ and $mw = O(q^{-\alpha\beta})$. Then, in view of the argument preceding (6.8), for any $\vec{\sigma} \in \{-1,1\}^m$ the regression function $\eta_{\vec{\sigma}}$ belongs to $\Sigma(\beta, L, \mathbf{R}^d)$ and Assumption (MA) is satisfied. We now check that $P_X$ obeys the strong density assumption. First, the density $\mu(x)$ equals a positive constant for $x$ belonging to the union of balls $\bigcup_{i=1}^m \mathcal{B}(z_i, (4q)^{-1})$, where $z_i$ is the center of $\mathcal{X}_i$ and $\mu(x) = (1 - mw)/(1 - mq^{-d}) = 1 + o(1)$, as $n \to \infty$, for $x \in A_0$. Thus, $\mu_{\min} \leq \mu(x) \leq \mu_{\max}$ for some positive $\mu_{\min}$ and $\mu_{\max}$. [Note that this construction does not allow one to choose any prescribed values of $\mu_{\min}$ and $\mu_{\max}$, because $\mu(x) = 1 + o(1)$. The problem can be fixed via a straightforward but cumbersome modification of the definition of $A_0$ that we skip here.] Second, the $(c_0, r_0)$-regularity of the support $A$ of $P_X$ with some $c_0 > 0$ and $r_0 > 0$ follows from the fact that, by construction, $\lambda(A \cap \mathcal{B}(x,r)) = (1 + o(1))\lambda([0,1]^d \cap \mathcal{B}(x,r))$ for all $x \in A$ and $r > 0$ (here again we skip the obvious generalization allowing to get any prescribed $c_0 > 0$). Thus, the strong density assumption is satisfied, and we conclude that $\mathcal{H} \subset \mathcal{P}_\Sigma$. Theorem 3.5 now follows from (6.8) if we choose $C'$ small enough.

Finally, we prove Theorem 4.1. Take $q = \lfloor Cn^{1/((2+\alpha)\beta+d)} \rfloor$, $w = C'q^{2\beta}/n$ and $m = q^d$ for some constants $C > 0$, $C' > 0$, and choose $A_0$ as a Euclidean ball contained in $\mathcal{X}_0$. As in the proof of Theorem 3.5, under the appropriate choice of $C$ and $C'$, the regression function $\eta_{\vec{\sigma}}$ belongs to $\Sigma(\beta, L, \mathbf{R}^d)$ and the margin Assumption (MA) is satisfied. Moreover, it is easy to see that the marginal distribution of $X$ obeys the mild density assumption [the $(c_0, r_0)$-regularity of the support of $P_X$ follows from considerations analogous to those in the proof of Theorem 3.5]. Thus, $\mathcal{H} \subset \mathcal{P}'_\Sigma$. Choosing $C'$ small enough and using (6.8), we obtain Theorem 4.1.

6.3. *Proof of Theorem* 4.2. We prove the theorem for $p < \infty$. The proof for $p = \infty$ is analogous. For any decision rule $f$ we set $d(f) \stackrel{\triangle}{=} R(f) - R(f^*)$ and

$$f^{**}(x, f) \stackrel{\triangle}{=} \begin{cases} f^*(x), & \text{if } \eta(x) \neq 1/2, \\ f(x), & \text{if } \eta(x) = 1/2, \end{cases} \qquad \forall x \in \mathbf{R}^d.$$

LEMMA 6.1. *Under Assumption* (MA), *for any decision rule $f$ we have*
(6.9) $$P_X(f(X) \neq f^{**}(X, f)) \leq Cd(f)^{\alpha/(1+\alpha)}.$$

PROOF. Note that $f^{**}(\cdot, f)$ is a Bayes rule, and following the same lines as in Proposition 1 of [20], we get $P_X(f(X) \neq f^{**}(X, f), \eta(X) \neq 1/2) \leq Cd(f)^{\alpha/(1+\alpha)}$. It remains to observe that $P_X(f(X) \neq f^{**}(X, f), \eta(X) \neq 1/2) = P_X(f(X) \neq f^{**}(X, f))$. □

For a Borel function $\bar{\eta}$ on $\mathbf{R}^d$ define $f_{\bar{\eta}} \triangleq \mathbb{1}_{\{\bar{\eta} \geq 1/2\}}$, $f_{\bar{\eta}}^*(\cdot) \triangleq f^{**}(\cdot, f_{\bar{\eta}})$ and
$$Z_n(f_{\bar{\eta}}) \triangleq [R_n(f_{\bar{\eta}}) - R_n(f_{\bar{\eta}}^*)] - [R(f_{\bar{\eta}}) - R(f_{\bar{\eta}}^*)] = [R_n(f_{\bar{\eta}}) - R_n(f_{\bar{\eta}}^*)] - d(f_{\bar{\eta}}).$$

Let $\eta_n$ be an element of $\mathcal{N}_{\varepsilon_n}$ such that $\|\eta_n - \eta\|_p \leq \varepsilon_n$, where $\|\cdot\|_p$ is the $L_p(\mathcal{C}, P_X)$ norm. It follows from the comparison inequality (5.3) that $d(f_{\eta_n}) \leq C\varepsilon_n^{(1+\alpha)p/(p+\alpha)} \triangleq \delta_n$. Set
$$\Delta_n = Cn^{-(1+\alpha)p/((2+\alpha)p+\rho(p+\alpha))}$$
(i.e., $\Delta_n$ is of the order of the desired rate). Fix $t > 0$ and introduce the set
$$\mathcal{N}_n^* = \{\bar{\eta} \in \mathcal{N}_{\varepsilon_n} : d(f_{\bar{\eta}}) \geq t\Delta_n\}.$$

For any $t > 0$ we have
$$\mathbf{P}(d(\hat{f}_n^s) \geq t\Delta_n)$$
$$\leq \mathbf{P}\left(\min_{\bar{\eta} \in \mathcal{N}_n^*}[R_n(f_{\bar{\eta}}) - R_n(f_{\eta_n})] \leq 0\right)$$
$$= \mathbf{P}\left(\min_{\bar{\eta} \in \mathcal{N}_n^*}[Z_n(f_{\bar{\eta}}) - Z_n(f_{\eta_n}) + d(f_{\bar{\eta}}) - d(f_{\eta_n})] \leq 0\right)$$
$$\leq \mathbf{P}\left(\min_{\bar{\eta} \in \mathcal{N}_n^*}[Z_n(f_{\bar{\eta}}) - Z_n(f_{\eta_n}) + d(f_{\bar{\eta}})/2 + t\Delta_n/2 - d(f_{\eta_n})] \leq 0\right)$$
$$\leq \mathbf{P}\left(\min_{\bar{\eta} \in \mathcal{N}_n^*}[Z_n(f_{\bar{\eta}}) + d(f_{\bar{\eta}})/2] \leq 0\right) + \mathbf{P}(Z_n(f_{\eta_n}) \geq t\Delta_n/2 - d(f_{\eta_n}))$$
$$\leq \mathbf{P}\left(\min_{\bar{\eta} \in \mathcal{N}_n^*}[Z_n(f_{\bar{\eta}}) + d(f_{\bar{\eta}})/2] \leq 0\right) + \mathbf{P}(Z_n(f_{\eta_n}) \geq t\Delta_n/2 - \delta_n).$$

Since $\Delta_n$ is of the same order as $\delta_n$, we can choose $t$ large enough to have $t\Delta_n/2 - \delta_n \geq t\Delta_n/4$. Thus,
$$\mathbf{P}(d(\hat{f}_n^s) \geq t\Delta_n) \leq \operatorname{card}\mathcal{N}_n^* \max_{\bar{\eta} \in \mathcal{N}_n^*} \mathbf{P}(Z_n(f_{\bar{\eta}}) \leq -d(f_{\bar{\eta}})/2)$$
$$+ \mathbf{P}(Z_n(f_{\eta_n}) \geq t\Delta_n/4)$$
$$\leq \exp(A'\varepsilon_n^{-\rho}) \max_{\bar{\eta} \in \mathcal{N}_n^*} \mathbf{P}(Z_n(f_{\bar{\eta}}) \leq -d(f_{\bar{\eta}})/2)$$
$$+ \mathbf{P}(Z_n(f_{\eta_n}) \geq t\Delta_n/4).$$

Note that for any Borel function $\bar{\eta}$ the value $Z_n(f_{\bar{\eta}})$ is an average of $n$ i.i.d. bounded and centered random variables whose variance does not exceed $P_X(f_{\bar{\eta}}(X) \neq f_{\bar{\eta}}^*(X)) \leq Cd(f_{\bar{\eta}})^{\alpha/(1+\alpha)}$ [cf. (6.9)]. Thus, using Bernstein's inequality we obtain

$$\mathbf{P}(-Z_n(f_{\bar{\eta}}) \geq a) \leq \exp\left(-\frac{Cna^2}{a + d(f_{\bar{\eta}})^{\alpha/(1+\alpha)}}\right) \qquad \forall a > 0.$$

Therefore, for $\bar{\eta} \in \mathcal{N}_n^*$,

$$\mathbf{P}(Z_n(f_{\bar{\eta}}) \leq -d(f_{\bar{\eta}})/2) \leq \exp(-Cnd(f_{\bar{\eta}})^{(2+\alpha)/(1+\alpha)})$$
$$\leq \exp(-Cn(t\Delta_n)^{(2+\alpha)/(1+\alpha)}).$$

Similarly, for $t > C$,

$$\mathbf{P}(Z_n(f_{\eta_n}) \geq t\Delta_n/4) \leq \exp\left(-\frac{Cn\Delta_n^2}{\Delta_n + d(f_{\eta_n})^{\alpha/(1+\alpha)}}\right)$$
$$\leq \exp\left(-\frac{Cn\Delta_n^2}{\Delta_n + \delta_n^{\alpha/(1+\alpha)}}\right)$$
$$\leq \exp\left(-Cn\Delta_n^{(2+\alpha)/(1+\alpha)}\right).$$

The result of the theorem follows now from the above inequalities and the relation $n\Delta_n^{(2+\alpha)/(1+\alpha)} \asymp \varepsilon_n^{-\rho}$.

Center for Education
 and Research in Informatics
Ecole Nationale des Ponts et Chaussées
19, rue Alfred Nobel
Cité Descartes, Champs-sur-Marne
77455 Marne-La-Valle
France
E-mail: audibert@certis.enpc.fr

Laboratoire de Probabilités et Modèles
 Aléatoires (UMR CNRS 7599)
Université Paris VI
4 pl. Jussieu, Boîte courrier 188
75252 Paris
France
E-mail: tsybakov@ccr.jussieu.fr